\documentclass[a4paper,11pt]{article}
\usepackage[french]{babel}
\usepackage{latexsym}
\usepackage{amssymb}
\usepackage{amsfonts}
\usepackage[ansinew]{inputenc}
\usepackage{mathrsfs}

\newtheorem{thm}{Théorème}[section]
\newtheorem{lem}{Lemme}[section]
\newtheorem{coro}{Corollaire}[section]
\newtheorem{propos}{Proposition}[section]
\newtheorem{defn}{Définition}[section]
\newtheorem{exe}{Exemple}[section]
\newtheorem{rmq}{Remarque}[section]

\newcommand{\pf}{{\bf Démonstration.~}} 

\newcommand{\thmref}[1]{Théorème~\ref{#1}}
\newcommand{\lemref}[1]{Lemme~\ref{#1}}
\newcommand{\corref}[1]{Corollaire~\ref{#1}}
\newcommand{\proposref}[1]{Proposition~\ref{#1}}
\newcommand{\defnref}[1]{Définition~\ref{#1}}

\input amssym.def
\input amssym
\input epsf

\font\ma=cmcsc10

\def\Proclaim #1. #2\par{\medbreak
  \noindent{\ma#1.\enspace}{\sl#2}\par\medbreak}

\def\scat#1#2%
 {\hbox{\vrule\vbox to#2pc{\hrule\hbox to#1pc{\hfil} 
  \vfil\hrule}\vrule}}

\def\Supp
 {\mathop{\rm Supp\,}\nolimits}

\def\conv
 {\mathop{\rm conv\,}\nolimits}

\def\Hom
 {\mathop{\rm Hom\,}\nolimits}

\def\card
 {\mathop{\rm card\,}\nolimits}

\def\codim
 {\mathop{\rm codim\,}\nolimits}

\def\re
 {\mathop{\rm Re\,}\nolimits}

\def\im
 {\mathop{\rm Im\,}\nolimits}

\def\aff
 {\mathop{\rm aff}\nolimits}

\def\Ch
 {\mathop{\rm Ch\,}\nolimits}

\def\hfl#1#2{\smash{\mathop{\hbox to 12mm{\rightarrowfill}}
 \limits^{\scriptstyle#1}_{\scriptstyle #2}}}

\def\vfl#1#2{\llap{$\scriptstyle #1$}\left\downarrow
 \vbox to 4mm{}\right.\rlap{$\scriptstyle #2$}}

\def\Hfl#1#2{\smash{\mathop{\hbox to 15mm{\rightarrowfill}}
 \limits^{\scriptstyle#1}_{\scriptstyle #2}}}

\begin{document}
\title{\bf Amibes de Sommes d'Exponentielles}
\author{\scshape{James Silipo}}
\date{7 Fevrier 2005}
\maketitle
\begin{abstract}
L'objectif de cet article est d'étudier la notion d'amibe au sens de Favorov pour les systèmes finis de sommes d'exponentielles à fréquences réelles et de montrer que, sous des hypothèses de généricité sur les fréquences, le complémentaire de l'amibe d'un système de~$(k+1)$ sommes d'exponentielles à fréquences réelles est un sous-ensemble $k$-convexe au sens d'Henriques.
\end{abstract}

\section{Introduction et énoncé du résultat principal.}

Soit~$P\subset\mathbb C[u_1^{\pm 1},\ldots,u_n^{\pm 1}]$ un système fini de polynômes de Laurent en $n$ variables
et~$V(P)$ son
ensemble de zéros dans le tore~$(\mathbb C^*)^n$; si~${\rm Log}$ est l'application de~$(\mathbb C^*)^n$
dans~$\mathbb R^n$ définie par
$$
{\rm Log}(u):=(\log\vert u_1\vert,\ldots,\log\vert u_n\vert)\,,\qquad u\in (\mathbb C^*)^n\,,
$$
l'{\sl amibe\/}~${\mathcal A}_P$  de~$P$ est l'image de~$V(P)$ par l'application~${\rm Log}$, soit
$$
{\mathcal A}_P:={\rm Log}\,V(P)\,.
$$
La notion d'amibe pour un seul polynôme de Laurent a été introduite par Gelfand,
Kapranov et Zelevinsky dans~\cite{GKZ} où l'on trouve exposées ses
propriétés fondamentales. Des études plus raffinées et des généralisations diverses de cette notion ont été faites par d'autres, parmi eux Forsberg, Passare, Rullg{\aa}rd, Tsikh (dont les travaux~\cite{For}, \cite{FPT}, \cite{Ru1}, \cite{Ru2}, \cite{PR} étudient
les relations entre l'amibe~${\mathcal A}_p$ d'un polynôme de Laurent~$p$, son polytope de
Newton~$\Gamma_p$ et les développements de Laurent de la fonction rationnelle~$1/p$) ou encore
Mikhalkin (qui donne dans \cite{Mi1} et \cite{Mi2} des applications et des généralisations de la notion d'amibe à la
géométrie des courbes réelles et tropicales). En parti\-culier, dans~\cite{GKZ} on trouve la preuve du fait suivant: 
\begin{propos}{\bf (\cite{GKZ})}\label{0-conv}
le complémentaire~${\mathcal A}_p^c$ de l'amibe d'un (seul) polynôme de Laurent~$p$ n'a q'un nombre fini de composantes connexes et chacune de ces composantes est convexe. 
\end{propos}
La \proposref{0-conv} cesse d'être vraie si l'on passe à un système~$P$ de polynômes de Laurent. En particulier les composantes connexes de~${\mathcal A}_P^c$ ne sont plus en général des ensembles convexes; cependant, Henriques~\cite{Hen} a observé que~${\mathcal A}_P^c$ vérifie une propriété plus faible qui s'exprime en des termes homologiques de la manière suivante. 
\begin{defn}{\bf (\cite{Hen})}\label{k-conv}
Soit~$k\in\mathbb N$,~$S\subseteq\mathbb R^n$ un~$(k+1)$-sous-espace affine orienté et~$Y\subseteq S$ un sous-ensemble. Une classe d'homologie (singulière) réduite dans~$\tilde H_k(Y,\mathbb Z)$ est dite {\rm non négative\/} si, pour tout point~$x\in S\setminus Y$, son image (sous le morphisme induit par l'inclusion) dans~$\tilde H_k(S\setminus\{x\},\mathbb Z)\simeq\mathbb Z$ est non négative. Le sous-ensemble des classes non négatives du groupe~$\tilde H_k(Y,\mathbb Z)$ est noté~$\tilde H^+_k(Y,\mathbb Z)$.
\par
Un sous-ensemble~$X\subseteq\mathbb R^n$ est dit~$k$-{\rm convexe} si pour tout~$(k+1)$-sous-espace affine orienté~$S\subset\mathbb R^n$, la classe nulle est la seule classe non négative de~$\tilde H_k(S\cap X,\mathbb Z)$ qui appartient au noyau du morphisme
$$
\tilde H_k(S\cap X,\mathbb Z)\rightarrow \tilde H_k(X,\mathbb Z)
$$ 
induit par l'inclusion.
\end{defn}
\begin{thm}{\bf (\cite{Hen})}\label{hen}
Soit~$P\subset\mathbb C[z^{\pm 1},\ldots,z^{\pm 1}]$ un système de polynômes de Lau\-rent tel que~$V(P)\subset(\mathbb C^*)^n$ a codimension~$(k+1)$. Alors,~${\mathcal A}_P^c$ est un sous-ensemble~$k$-convexe.
\end{thm}
Cet énoncé peut se lire comme un résultat d'injectivité partielle du morphisme
$$
\iota_{k,S}:\tilde H_k(S\cap {\mathcal A}_P^c,\mathbb Z)\rightarrow \tilde H_k({\mathcal A}_P^c,\mathbb Z)
$$
pour chaque~$(k+1)$-sous-espace affine orienté~$S\subset\mathbb R^n$. Si~$k=0$, les morphismes~$\iota_{0,S}$ correspondant sont effectivement tous injectifs (et dans ce cas le \thmref{hen} se réduit à la \proposref{0-conv}), par contre, dès que~$k>0$, les morphismes~$\iota_{k,S}$ ne le sont plus que dans un sens conjectural, (voir~\cite{Mi2}). 
\par
Les travaux de Ronkin et Favorov autour des amibes soulèvent des questions nouvelles et tout particulièrement intéressantes dans l'étude des certains sous-ensembles analytiques globaux de~$\mathbb C^n$. En fait, les articles \cite{Ron} et \cite{Fav} adaptent la notion d'amibe au cadre des fonctions holomorphes presque périodiques définies dans les domaines de~$\mathbb C^n$ du type~$T_\Omega:=\mathbb R^n+i\Omega\,$,~$\Omega$ étant un ouvert de~$\mathbb R^n$. Il s'agit de la classe~
$
AP(T_\Omega)
$
des fonctions~$g\in{\mathcal O}(T_\Omega)$ telles que l'ensemble~$\{g(z+t)\in{\mathcal O}(T_\Omega)\mid t\in\mathbb R^n\}$ est relativement compacte dans la topologie~$\tau(T_\Omega)$ induite sur~${\mathcal O}(T_\Omega)$ par la convergence uniforme sur les sous-domaines du type~$T_D$, avec~$D\Subset\Omega$. 
\begin{defn}\label{fav}
Soit~$\Omega\subset\mathbb R^n$ un ouvert non vide. L'amibe d'un système fini~$G\subset AP(T_\Omega)$ est le sous-en\-semble de~$\mathbb R^n$ donné par
$$
A_G:=\overline{\im V(G)}\,,
$$
o\`u~$V(G)$ dénote l'ensemble de zéros de~$G$ dans~$T_\Omega$ et~$\im:T_\Omega\rightarrow\Omega$ est l'application de prise de partie imaginaire sur chaque coordonnée.
\end{defn}
Dans Favorov~\cite{Fav} on trouve la~\defnref{fav} dans le cas d'un système réduit à une seule fonction et afin d'éviter toute ambiguïté entre les notations~${\mathcal A}_{P}$ et~$A_{G}$, on va dorénavant indiquer les amibes au sens de Favorov (soit au sens de la \defnref{fav}) par le symbole~${\mathcal F}_{G}$.
\par
Un cas bien particulier (mais néanmoins très important\footnote{Un résultat profond de la théorie des fonctions holomorphes presque périodiques (le théorème d'approximation de Bochner-Fejér) assure que toute fonction~$g\in AP(T_\Omega)$ est la limite dans la topologie~$\tau(T_\Omega)$ d'une suite convergeante de sommes d'exponentielles à fréquences imaginaires pures.}) de fonctions de~$AP(T_{\mathbb R^n})=AP(\mathbb C^n)$ est celui des sommes d'exponentielles à fréquences imaginaires pures, soit les fonctions de la forme
\begin{eqnarray}
g(z)
=
\sum_{\lambda\in\Lambda}c_\lambda \,e^{i\langle z,\lambda\rangle}
=
\sum_{\lambda\in\Lambda}c_\lambda \,e^{\langle z,-i\lambda\rangle}
\end{eqnarray}
o\`u~$z\in\mathbb C^n$,~$\Lambda\subset\mathbb R^n$ est un ensemble fini et~$c_\lambda\in\mathbb C^*$ pour tout~$\lambda\in\Lambda$ (les vecteurs~$-i\lambda\in i\mathbb R^n$ étant les {\sl fréquences\/} de~$g$). Toutefois, dans cet article on va plutôt travailler avec les systèmes finis de sommes d'exponentielles à fréquences réelles, soit les systèmes finis de fonctions du type
$$
f(z)=g(-iz)\,,
$$
o\`u~$g$ est de la forme~$(1)$ ci-dessus, donc pour un tel système on va dorénavant assumer la définition suivante.
\begin{defn} Soit~$F$ un système de sommes d'exponentielles à fréquences réelles. L'amibe au sens de Favorov de~$F$ est l'ensemble
$$
{\mathcal F}_F:=\overline{\re V(F)}\,,
$$
où~$V(F)$ dénote l'ensemble de zéros de~$F$ dans~$\mathbb C^n$ et~$\re:\mathbb C^n\rightarrow\mathbb R^n$ est l'application de prise de partie réelle sur chaque coordonnée.
\end{defn}
\par
Comme remarqué dans~\cite{Ron} ou~\cite{Fav}, si~$g\in AP(T_\Omega)$, puisqu'elle est holomorphe, chaque composante connexe de l'ensemble~${\mathcal F}_g^c\cap\Omega$ est aussi convexe. En outre, si~$g$ (resp.~$f$) est une somme d'exponentielles à fréquences imaginaires pures (resp. réelles), l'ensemble~$\mathbb R^n\setminus\overline{\im V(g)}$, (resp.~$\mathbb R^n\setminus\overline{\re V(f)}$) n'a qu'un nombre fini de composantes connexes convexes, donc la \proposref{0-conv} se traduit mot à mot au cadre des amibes des sommes d'exponentielles à fréquences imaginaires pures (resp. réelles). 
\par
Si l'on passe aux systèmes finis des fonctions de~$AP(T_\Omega)$, la structure des amibes devient considérablement plus compliquée. Cependant, dans le cadre des systèmes finis de sommes d'exponentielles à fréquences réelles (resp. imaginaires pures), la théorie développée par Ka\-zar\-nov\-ski{\v\i}~\cite{Ka1} permet, d'une part, de mieux comprendre la structure des amibes au sens de Favorov associées à ces systèmes et, d'autre part, d'adapter au même cadre le résultat d'Henriques~\cite{Hen}. Pour énoncer notre résultat on a besoin de quelques notations qui seront détaillées dans les sections suivantes.
\par
Si~$F$ est un système fini de sommes d'exponentielles à fréquences réelles, on associe à~$F$ une famille~$\{F_\chi\}_\chi$ de systèmes ``perturbés'' du système~$F$, l'indice~$\chi$ parcourant un certain groupe de caractères associé à~$F$.
On introduit ainsi une nouvelle notion d'{\sl amibe\/} en posant
$$
{\mathcal Y}_F
:=
\bigcup_{\chi}
\,\re V(F_\chi)\,
$$
et l'on obtient le résultat suivant (voir Théorème~3.1 et Théorème~5.1 respectivement).
\par\bigskip\noindent
{\bf Résultat.} {\sl Soit~$F$ un système constitué par~$(k+1)$ sommes d'exponentielles à fré\-quences réelles génériques, alors}
\begin{description}
\item [$(a)$]
$$
{\mathcal Y}_F
=
\mathbb R^n\cap\bigcup_{\chi} V(F_\chi)
=
\overline{\re V(F)}\,,
$$
\noindent
{\sl en particulier l'amibe~${\mathcal Y}_{F}$ coïncide avec l'amibe~${\mathcal F}_{F}$ au sens de Favorov;}
\end{description}
\begin{description}
\item[$(b)$] 
{\sl ~~le complé\-mentaire~${\mathcal F}^c_F$ de l'amibe de~$F$ est~$k$-convexe dans~$\mathbb R^n$.}
\end{description}
\par
\bigskip
La partie~$(a)$ fournit un expression plus concrète de l'adhérence de l'ensemble~$\re V(F)$ et elle implique, entre autres, que 
$$
\overline{\re V(F)}
=
\overline{\re V(F_\chi)}\,,
$$
pour tout~$\chi\,$. La partie~$(b)$ constitue le pendant du \thmref{hen} dans le cadre exponentiel. Les preuves de~$(a)$ et de~$(b)$ utilisent la technique de perturbation par caractères introduite depuis longtemps par A.~Yger dans les travaux~\cite{Y1} et~\cite{Y2} (puis utilisés par C.~Berenstein et A.~Yger) pour montrer que certains systèmes d'équations de convolution possédaient la propriété de la synthèse spectrale. En ce sens, la présentation de l'amibe de Favorov donnée dans~$(a)$ pourrait s'avérer intéressante du point de vue des questions de petits dénominateurs inhé\-rentes aux systèmes à fréquences réelles non commensurables.

\section{Sommes d'exponentielles: définitions et notations}
Dans cette section on rappelle toutes les notions et tous les résultats autour des sommes d'exponentielles utiles dans la suite.
\par
Soit~$n\in \mathbb N^*$ fixé et~${\mathcal O}(\mathbb C^n)$ la~$\mathbb C$-algèbre des fonctions holomorphes sur~$\mathbb C^n$.
Une {\sl som\-me d'expo\-nentiel\-les\/} sur~$\mathbb C^n$ est
un élément de la sous-algèbre~${\mathcal S}_n$ de~${\mathcal O}(\mathbb C^n)$ engendrée, en tant que sous-espace vectoriel complexe, par
les fonctions de la forme~$e^{\langle z,\lambda\rangle}$, o\`u~$\lambda\in\mathbb C^n$.
${\mathcal S}_n^*$ dénote l'ensemble des sommes d'exponentielles non nulles.
\par
Si~$f\in{\mathcal S}_n^*$, le {\sl spectre\/} de~$f$ est le plus petit
sous-ensemble~$\Lambda_f$ de~$\mathbb C^n$ tel que~$f$ appartient au sous-espace vectoriel
de~${\mathcal S}_n$ engendré par l'ensemble des monômes exponentiels~$\{ e^{\langle z,\lambda\rangle}\mid \lambda\in\Lambda_f\}$, (il s'agit d'un ensemble bien défini puisque la famille~$\{e^{\langle z,\lambda\rangle}\}_{\lambda\in{\mathbb C}^n}$ est une base de~${\mathcal S}_n$ sur~$\mathbb C$),
les {\sl fréquences\/} de~$f$ sont les éléments de son spectre~$\Lambda_f$.
Le {\sl polytope de Newton\/} de~$f\in{\mathcal S}_n^*$ est l'enveloppe convexe
$$
\Gamma_f:=\conv \Lambda_f\subset\mathbb C^n
$$
de son spectre~$\Lambda_f$. À toute~$f\in{\mathcal S}_n^*$ on associe la fonction réelle~$k_f$ donnée, pour~$z\in\mathbb C^n$, par
$$
k_f(z)
:=
\sup_{\lambda\in\Lambda_f} e^{\re\langle z,\lambda\rangle}\,;
$$
la fonction~$k_f$ n'est rien d'autre que l'exponentielle de la fonction de support du polytope de Newton de~$f$, calculée par rapport au produit scalaire~$\re\langle\,,\rangle$ sur~$\mathbb C^n$.
\par
Dans cet article on utilisera des sous-algèbres de~${\mathcal S}_n$, à savoir les sous-algèbres du type~${\mathcal S}_{n,\mathbb G}$
constituées par les sommes d'exponentielles à fré\-quences dans un sous-groupe additif~$\mathbb G$ de~$\mathbb C^n$, souvent~$\mathbb G$ sera~$\mathbb Z^n,\mathbb Q^n,\mathbb R^n$ ou~$i\mathbb R^n$, on parlera ainsi de sommes d'exponentielles {\sl à fréquences entières, rationnelles, réelles ou imaginaires pures\/}; on note que pour tout~$\mathbb G$, on a~${\mathcal S}_{n,\mathbb G}^*={\mathcal S}_n^*\cap{\mathcal S}_{n,\mathbb G}$.
\par
Un {\sl système de sommes d'exponentielles\/}~(en abrégé SSE) est un
sous-en\-semble non vide et fini~$F$ de~${\mathcal S}_n^*$. Pour un tel système~$F$
on pose
$$
\Gamma_F:=\sum_{f\in F}\Gamma_f
$$
(la somme au deuxième membre étant prise au sens de Minkowski).
L' {\sl ensemble des spectres\/} de~$F$ est l'ensemble~$\{\Lambda_f\mid f\in F\}$ et les
{\sl fréquences\/} de~$F$ sont les éléments de l'union des spectres des~$f\in F$;~$F$ est dit à {\sl fréquences entières, rationnelles, réelles ou imaginaires pures} si chaque~$f\in F$ l'est.
On note, re\-specti\-ve\-ment,~$\Xi_F$,~${\rm vect}_\mathbb Q \Xi_F$ 
et~${\rm vect}_{\mathbb R}\Xi_F$ le sous-groupe additif, le~$\mathbb Q$-sous-espace vectoriel et le~$\mathbb R$-sous-espace vectoriel de~$\mathbb C^n$ engendrés
par les fré\-quences de~$F$. 
\par
Si~$\mathbb G\subset \mathbb C^n$ est un sous-groupe additif qui contient les fréquences de~$F$, pour tout homomorphisme~$\chi$ de groupes abéliens, du groupe additif~$\mathbb G$ à valeurs dans le groupe multiplicatif~$\mathbb S^1$ des nombres complexes de module égale à~$1$, ~$\chi\in\Ch\mathbb G:=\Hom_{\mathbb Z}(\mathbb G,\mathbb S^1)$, on introduit le SSE
$$
F_\chi:=\{f_\chi\in{\mathcal S}_n^*\mid f\in F\}\,, 
$$
où, pour tout~$f\in F$, on a posé
$$
f_\chi(z):=\sum_{\lambda\in\Lambda_f}c_\lambda \chi(\lambda)e^{\langle z,\lambda\rangle}\,.
$$ 
On observe que le groupe abélien~$\mathbb S^1$ est divisible, donc il est un objet injectif dans la catégorie des groupes abéliens (voir \cite{Eis}), ce qui est équivalent à la surjectivité de l'homomorphisme de
restriction
$$
\rho:\Ch \mathbb G \longrightarrow\Ch\Xi_f\,.
$$
On en déduit que l'ensemble~
$
\{F_\chi\subset{\mathcal S}_n^*\mid\chi\in\Ch\mathbb G\}
$
ne dépende que de~$F$, quel que soit le sous-groupe additif~$\mathbb G$ de~$\mathbb C^n$ contenant les fréquences de~$F$. 
\par
\`A tout SSE~$F$ on associe la fonction réelle bornée~$K[F]$ donnée, pour tout~$z\in\mathbb C^n$, par
$$
K[F](z)
:=
\sum_{f\in F}
{\vert f(z)\vert
\over
k_f(z)}\,
$$
et à toute face~$\Delta=\sum_{f\in F}\Delta_f$ de~$\Gamma_F$ on associe le SSE 
$$
F^\Delta:=\{f^\Delta\in{\mathcal S}_n^*\mid f\in F\}\,,
$$
appelé {\sl $\Delta$-trace\/} de~$F$, obtenu en posant, pour~$f\in F$,
$$
f^\Delta(z):=
\sum_{\lambda\in \Lambda_f\cap \Delta_f} c_\lambda e^{\langle z,\lambda\rangle}
\,.
$$
\par
Si~$\Delta$ est une face d'un polytope~$\Gamma\subset\mathbb C^n$, on note~$\aff_{\mathbb C}\Delta$ le sous-espace affine complexe de~$\mathbb C^n$ engendré par~$\Delta$. 
\par
Les notations que l'on vient de préciser, permettent de reprendre certaines notions introduites par Kazarnovski{\v\i} \cite{Ka1}. 
\begin{defn}{\bf (\cite{Ka1})}
Un SSE~$F$ est dit {\rm régulier\/} s'il existe~$\varepsilon>0$ tel que,
pour chaque~$\Delta\preccurlyeq\Gamma_F$ avec~$\dim_{\mathbb C}(\aff_{\mathbb C}\Delta)<\card F$, on a~$K[F^\Delta]\geqslant\varepsilon$.
\end{defn}
\par
\begin{thm}{\bf (\cite{Ka1})}\label{kaza}
Soit~$F$ un SSE régulier, alors l'ensemble~$V(F)$, des zéros de~$F$ dans~$\mathbb C^n$, est non vide si et seulement si~$\dim_{\mathbb C}(\aff_{\mathbb C}\Gamma_F)\geqslant\card F$ et dans ce cas sa codimension est égale à~$\card F$.\hfill$\square$
\end{thm}
\par 
\begin{defn}{\bf (\cite{Ka1})}\label{sf}
L'ensemble des spectres d'un~SSE~$F$ est dit {\rm fermé} si pour toute face~$\Delta=\sum_{f\in F}\Delta_f$ de~$\Gamma_F$ telle que~$\dim_{\mathbb C}(\aff_{\mathbb C}\Delta)<\card F$, il existe~$f$ dans~$F$ pour
lequel~$\Delta_f$ soit réduit à un point.
\end{defn}
\par
\begin{propos}\label{kaza2}{\bf (\cite{Ka1})}
Un SSE dont l'ensemble des spectres est fermé est un SSE régulier.\hfill$\square$
\end{propos}
La condition de la~\defnref{sf} ne regarde que les spectres du système (ou m\^eme juste les sommets des polytopes~$\Gamma_f$,~$f\in F$), donc si l'on fixe un nombre~$r\leqslant n$ de spectres dans~$\mathbb C^n$ (resp. dans~$\mathbb R^n$) ainsi qu'un~$r$-uplet~$(\ell_1,\ldots,\ell_r)\in\mathbb N^r$, on peut montrer que l'ensemble des spectres d'un SSE constitué par~$r$ sommes d'exponentielles dont les spectres comportent respectivement~$\ell_1,\ldots,\ell_r$ fréquences, est génériquement fermé. On en déduit que, si~$F$ est un SSE dont l'ensemble des spectres est fermé, alors, pour tout~$\chi\in\Ch\Xi_F$, il en est de m\^eme ainsi de~$F_\chi$; en particulier, si~$V(F)\neq\varnothing$
(soit si~$\dim_{\mathbb C}(\aff_{\mathbb C}\Gamma_F)\geqslant\card F$) alors~$V(F_\chi)\neq\varnothing$ et
$$
\codim V(F)=\card F=\card F_\chi=\codim V(F_\chi)\,,
$$
pour tout~$\chi\in\Ch\Xi_F$.
\par
\section{Amibes: définition et premières propriétés.}
Suite à une idée d'Alain Yger \cite{Y1},~\cite{Y2}, on propose ici une nouvelle définition d'amibe pour les systèmes finis de sommes d'exponentielles à fré\-quences réelles.\footnote{Cette notion d'amibe pourrait s'adapter au cas plus général des systèmes finis de fonctions holomorphes presque périodiques dans les domaines tubulaires de~$\mathbb C^n$.} On verra en suite sous quelles conditions cette notion d'amibe coïncide avec celle due à Favorov.
\par
\begin{defn}\label{amibe}
Soit~$F$ un {\rm SSE}
à fréquences réelles,~$\mathbb G\subset\mathbb C^n$ un sous-groupe contenant les fréquences de~$F$.
On appelle {\rm amibe} de~$F$ le sous-ensemble~${\mathcal Y}_F$ de~$\mathbb R^n$ défini par
$$
{\mathcal Y}_F:=\bigcup_{\chi\in\Ch \mathbb G}\re V(F_\chi)\,.
$$
\end{defn}
\par\noindent
L'amibe~${\mathcal Y}_F$ est bien définie car l'ensemble~$\{F_\chi\subset{\mathcal S}_{n,\mathbb R^n}^*\mid \chi\in\Ch\mathbb G\}$ utilisé dans sa définition est indépendant du choix du sous-groupe additif~$\mathbb G\subset\mathbb C^n$ parmi ceux qui contiennent les fréquences de~$F$. Ceci nous autorise entre autre à représenter l'amibe~${\mathcal Y}_F$ à l'aide du groupe~$\mathbb G$ qui nous convient le plus. 
On remarque aussi que si~$\chi\in\mathbb G$ alors~$F$ et~$F_\chi$ ont les mêmes fréquences et donc les mêmes amibes:~${\mathcal Y}_F={\mathcal Y}_{F_\chi}$. 
\par
\begin{propos}\label{propr} Soit~$F$ un SSE à fréquences réelles, alors
\item[$(i)$] $${\mathcal Y}_F=\mathbb R^n\cap\bigcup_{\chi\in\Ch\Xi_F}V(F_\chi)\,,$$
\item[$(ii)$] ${\mathcal Y}_F$ est un sous-ensemble fermé dans~$\mathbb R^n$.
\end{propos}
\pf
$(i)$ Soit~$r$ le rang\footnote{Puisque le groupe additif~$\mathbb R^n$ n'a pas de torsion, ses sous-groupes de type fini sont libres.} de~$\Xi_F$ et~$\{\omega_1,\ldots,\omega_r\}$ un système de générateurs libres
de~$\Xi_F$, alors, pour tout~$f\in F$, on a l'expression
$$
f(z)
=
\sum_{k\in A_f}
a_{f,k}
\Big(e^{i\langle\omega_1,\im z\rangle}\Big)^{k_1}
\!\!\!\!\!
\cdots
\Big(e^{i\langle\omega_r,\im z\rangle}\Big)^{k_r}
e^{\langle k_1\omega_1+\cdots+k_r\omega_r,\re z\rangle}
$$
où~$A_f\subset\mathbb Z^r$ est un sous-ensemble fini et~$a_{f,k}\in\mathbb C^*$ pour tout~$k\in A_f$.
Si~$\xi\in{\mathcal Y}_F$, il existe un~$\chi\in\Ch\Xi_F$ et un~$\eta\in\mathbb R^n$ tels que~$f_\chi(\xi+i\eta)=0$, pour tout~$f\in F$. Par conséquent, si pour tout~$1\leqslant j\leqslant r$,~$\theta_j$ désigne la détermination principale de l'argument de~$\chi(\omega_j)$ on a
$$
f_\chi(\xi+i\eta)
=
\sum_{k\in A_f}
a_{f,k}
\Big(e^{i(\theta_1+\langle\omega_1,\eta\rangle)}\Big)^{k_1}
\!\!\!\!\!
\cdots
\Big(e^{i(\theta_r+\langle\omega_r,\eta\rangle)}\Big)^{k_r}
e^{\langle k_1\omega_1+\cdots+k_r\omega_r,\xi\rangle}
$$
pour tout~$f\in F$ donc, si~$\chi^\prime$ dénote le caractère de~$\Xi_F$ donné, pour~$1\leqslant j\leqslant r$, par~
$
\chi^\prime(\omega_j)=e^{i\langle \eta,\omega_j\rangle}\,,
$
on aura~$f_{\chi\chi^\prime}(\xi)=f_\chi(\xi+i\eta)=0$, pour tout~$f\in F$, soit~$\xi\in V(F_{\chi\chi^\prime})$. L'autre inclusion est triviale. 
\par
$(ii)$ Si $(\xi_q)_{q\in\mathbb N}\subset{\mathcal Y}_F$ est une suite convergeante 
vers~$\xi\in\mathbb R^n$, pour tout indice~$q\in\mathbb N$, il existe, gr\^ace à~$(i)$, un~$\chi_q\in\Ch \Xi_F$ tel que,~$f_{\chi_q}(\xi_q)=0\,,$ pour tout~$f\in F$. 
En vertu de la compacité de~$\Ch\Xi_F$, la suite~$(\chi_q)_{q\in\mathbb N}\subset\Ch \Xi_F$ admet une 
sous-suite~$(\tilde\chi_{q_m})_{m\in\mathbb N}$ qui converge 
vers un caractère~$\chi\in\Ch \Xi_F$, donc 
$$
f_\chi(\xi)=\lim_{m\to \infty}f_{\chi_{q_m}}(\xi_{q_m})=0\,,
$$
pour tout~$f\in F$, soit~$\xi\in{\mathcal Y}_F$.\hfill$\square$
\goodbreak
\bigskip\noindent
\begin{rmq}
{\rm 
Si~$F$ est un SSE à fréquences réelles, on a évidemment
$$
{\mathcal Y}_F
=
\bigcup_{\chi\in\Ch\Xi_F}
{\mathcal F}_{F_{_\chi}}
\qquad
{\rm et}
\qquad
{\mathcal Y}_F^c
=
\bigcap_{\chi\in\Ch\Xi_F}
{\mathcal F}^c_{F_{_\chi}}\,,
$$ 
ainsi que les inclusions~
$
\re V(F_\chi)\subseteq {\mathcal F}_{F_{_\chi}}\subseteq {\mathcal Y}_F\,,
$
en général strictes et valables pour tout~$\chi\in\Ch\Xi_F$.}
\end{rmq}
On va maintenant s'intéresser plus en détail aux rapports entre la notion d'amibe que l'on vient de définir et celle due à Favorov. Le~\thmref{proprF} fait le lien entre les deux notions mais sa preuve utilise une version multidimensionnelle d'un théo\-rème dit {\sl d'Approximation de Kronecker\/} (déjà utilisée par Ronkin~\cite{Ron}) dont on préfère ajouter ici une démonstration. On montre d'abord un lemme.
\begin{lem}\label{lemmey}
Soit~$H\subset\mathbb R^r$ un sous-groupe (additif) fermé. Alors~$H=\mathbb R^r$ ou bien il existe une forme $\mathbb R $-linéaire~$\psi\not\equiv0$ sur~$\mathbb R^r$ telle que~$\psi(H)\subseteq\mathbb Z$. 
\end{lem}
\pf Si~$r=1$ le lemme est une simple consé\-quence du fait bien connu qu'un sous-groupe additif de~$\mathbb R$ est soit dense soit discret. Par récurrence, on suppose que le lemme soit vrai dans~$\mathbb R^s$, pour tout~$s<r$. Or, si~$H\neq\mathbb R^r$, il existe un sous-espace linéaire~$S\subset\mathbb R^r$, avec~$0<\dim S<r$, tel que l'image de~$H$, sous la projection orthogonale~$\pi_S:\mathbb R^r\to S$, constitue un sous-groupe discret de~$S$. Puisque~$\pi_S(H)$ est discret dans~$S$, il existe par hypothèse de récurrence une forme~$\mathbb R$-linéaire~$\psi_S\not\equiv 0$ sur~$S$ telle que~$\psi_S(\pi_S(H))\subseteq \mathbb Z$. La forme~$\mathbb R$-linéaire~$\psi:=\psi_S\circ\pi_S$ vérifie le lemme pour~$H\subset\mathbb R^r$.\hfill$\square$
\par
\begin{thm}\label{multikro}
Soient $\omega_1,\ldots,\omega_r\in\mathbb R^n$ vecteurs linéairement indépendants sur $\mathbb Z$. Alors le sous-groupe additif
$$
G
:=
\{x\in\mathbb R^r\mid x_j=\langle t,\omega_j\rangle+p_j\quad\hbox{\rm où}\quad t\in\mathbb R^n,\;p_j\in\mathbb Z\quad{\rm et}\quad j=1,\ldots,r\}
$$
est dense~\footnote{On observe que des conditions diophantiennes portant sur les~$\omega_1,\ldots,\omega_r$ ``freinent'' la vitesse de l'approximation du point courant de~$\mathbb R^r$ par des points de~$G$.}dans~$\mathbb R^r$.
\end{thm}
\pf
Soit~$H$ l'adhérence de~$G$.~$H$ est aussi un sous-groupe et on va montrer que~$H=\mathbb R^m$. Si, par l'absurde,~$H\subsetneq\mathbb R^r$, le~\lemref{lemmey} implique qu'il existe une forme $\mathbb R$-linéaire~$\psi\not\equiv 0$ sur~$\mathbb R^r$ telle que~$\psi(x)\in\mathbb Z$ pour tout~$x\in H$ et donc {\it à fortiori} pour tout~$x\in G$. En évaluant~$\psi$ sur un point~$p$ de~$\mathbb Z^r\subset G$ on trouve
$$
\psi(p)=q_1p_1+\cdots+q_rp_r\,,
$$ 
pour certains~$q_1,\ldots,q_r\in\mathbb Z$ non tous nuls, d'autre part en l'évaluant sur un point du type
$$
x_t
=
(\langle t,\omega_1\rangle,\ldots,\langle t,\omega_r\rangle)\in G\,,
$$
où~$t\in\mathbb R^n$ est arbitraire, on obtient 
$$
\psi(x_t)
=
q_1\langle t,\omega_1\rangle+\cdots+q_r\langle t,\omega_r\rangle
=
\langle t,q_1\omega_1+\cdots+q_r\omega_r\rangle\in\mathbb Z
$$
pour tout~$t\in\mathbb R^n$, ce qui est possible si et seulement si
$$
q_1\omega_1+\cdots+q_r\omega_r=0\,,
$$
d'où la contradiction.\hfill$\square$
\par
\bigskip
Les preuves du~\lemref{lemmey} et du~\thmref{multikro} ont été obtenues en modifiant celles que l'on trouve dans~\cite{Mey}. 
\goodbreak
\par
\begin{thm}\label{proprF} 
Soit~$F=\{f_1,\ldots,f_s\}\subset{\mathcal S}^*_{n,\mathbb R^n}\,$,~$s\leqslant n$, tel que~$V(F_\chi)=\varnothing$ pour tout~$\chi\in\Ch\Xi_F$, ou bien tel que~$\codim V(F_\chi)=s$ pour tout~$\chi\in\Ch\Xi_F$.
Alors
$$
{\mathcal Y}_F=\overline{\re V(F)}={\mathcal F}_F\,.
$$
\end{thm}
\pf Si~$V(F_\chi)$ est vide pour tout~$\chi\in\Ch\Xi_F$, le théorème est vrai trivialement. 
Soit donc~$\codim V(F_\chi)=s$, (en particulier~$V(F_\chi)\neq\varnothing$), pour tout~$\chi\in\Ch\Xi_F$. L'inclusion~${\mathcal Y}_F\supseteq \overline{\re V(F)}$ est évidente, on doit donc montrer que~$\re V(F_\chi)\subseteq\overline{\re V(F)}$, pour tout~$\chi\in\Ch\Xi_F$. 
On suppose, par l'absurde, qu'il existe un~$\chi\in\Ch\Xi_F$ et un point~$x\in{\mathbb R}^n\cap V(F_\chi)$ qui n'est pas adhérent à~$\re V(F)$. Cela signifie qu'il existe un~$\varepsilon>0$ tel et que la bande
$$
D:=\{z\in\mathbb C^n\mid\,\Vert x-\re z\Vert<\varepsilon\}
$$
ne contient pas de zéros de~$F$. Si~$\{\omega_1,\ldots,\omega_r\}$ est un système de générateurs libres du groupe~$\Xi_F$ et, pour~$1\leqslant\ell\leqslant r$,~$\theta_\ell$ est la détermination principale de l'argument de~$\chi(\omega_\ell)$, on a les expressions suivantes pour tout~$1\leqslant j\leqslant s$,
$$
f_j(z)
=
\sum_{k\in A_j}
a_{j,k}
\Big(e^{i\langle\omega_1,\im z\rangle}\Big)^{k_1}
\!\!\!\!\!
\cdots
\Big(e^{i\langle\omega_r,\im z\rangle}\Big)^{k_r}
e^{\langle k_1\omega_1+\cdots+k_r\omega_r,\re z\rangle}
$$
et
$$
f_{j,\chi}(z)
=
\sum_{k\in A_j}
a_{j,k}
\Big(e^{i(\theta_1+\langle\omega_1,\im z\rangle)}\Big)^{k_1}
\!\!\!\!\!
\cdots
\Big(e^{i(\theta_r+\langle\omega_r,\im z\rangle)}\Big)^{k_r}
e^{\langle k_1\omega_1+\cdots+k_r\omega_r,\re z\rangle},
$$
où~$A_j$ est un sous-ensemble fini de~$\mathbb Z^r$ et~$a_{j,k}\in\mathbb C^*$ pour tout~$k\in A_j$.
Le~\thmref{multikro} implique qu'il existe une suite~$(t_m)\subset\mathbb R^n$ telle que, pour tout~$1\leqslant\ell\leqslant r$, on ait
$$
\lim_{m\to+\infty}
\langle\omega_\ell,t_m\rangle
=
\theta_\ell
\qquad
{\rm mod.}\;2\pi\mathbb Z\,,
$$
ce qui fait que, pour tout~$1\leqslant j\leqslant s$ et tout~$z\in\mathbb C^n$,
$$
\lim_{m\to+\infty}f_j(z+it_m)=f_{j,\chi}(z)\,.
$$
Pour tout~$m\in \mathbb N$, tout~$1\leqslant j\leqslant s$ et tout~$z\in\mathbb C^n$, on pose
$$
g_{j,m}(z)
:=
f_j(z+it_m)\,,
$$
donc~
$
\lim_{m\to+\infty}
g_{j,m}
=
f_{j,\chi}\,,
$
ce qui fait que le système~$G_m=\{g_{1,m},\ldots,g_{s,m}\}$ ``tend\footnote{L'idée d'approcher les~$f_{j,\chi}$ par des translatées des~$f_j$ à l'aide du~\thmref{multikro} a été déjà exploitée par Ronkin dans~\cite{Ron}.}'' vers~$F_\chi$ si~$m$ tend vers l'infini.
Or, puisque~$V(F)\cap D=\varnothing$, on a aussi~$V(G_m)\cap D=\varnothing$, pour tout~$m\in\mathbb N$, mais comme par construction~$\codim V(G_m)=s=\codim V(F_\chi)$ pour tout~$m\in\mathbb N$, la version plusieurs variables du théorème de Rouché fait que~$V(F_\chi)\cap D=\varnothing$ aussi. Ceci est absurde car par hypothèse~$x\in V(F_\chi)\cap D$.  \hfill$\square$
\begin{coro}
Soit~$F=\{f_1,\ldots,f_s\}\subset{\mathcal S}^*_{n,\mathbb R^n}\,$,~$s\leqslant n$, tel que~$V(F_\chi)=\varnothing$ pour tout~$\chi\in\Ch\Xi_F$, ou bien tel que~$\codim V(F_\chi)=s$ pour tout~$\chi\in\Ch\Xi_F$.
Alors
$$
\overline{\re V(F)}
=
\overline{\re V(F_\chi)}\,,
$$
pour tout~$\chi\in\Ch\Xi_F$.
\end{coro}
\pf Pour tout~$\chi\in\Ch\Xi_F$ le système~$F_\chi$ vérifie les mêmes hypothèses que~$F$ donc le \thmref{proprF} fait qu'on ait aussi
$$
{\mathcal Y}_{F_\chi}=\overline{\re V(F_\chi)}={\mathcal F}_{F_\chi}\,,
$$
d'autre part~${\mathcal Y}_F={\mathcal Y}_{F_\chi}$, donc~
$
\overline{\re V(F)}
=
\overline{\re V(F_\chi)}\,.
$\hfill$\square$
\begin{coro}
Soit~$F\subset{\mathcal S}^*_{n,\mathbb R^n}$ un SSE dont l'ensemble des spectres est fermé, alors
$$
{\mathcal Y}_F=\overline{\re V(F)}={\mathcal F}_F\,.
$$
\end{coro}
\pf Il suffit de remarquer qu'à cause de la~\proposref{kaza2}, pour tout caractère~$\chi\in\Ch\Xi_F$, le système~$F_\chi$ est régulier. Or pour le \thmref{kaza} on n'a plus que deux possibilités, ou bien~$V(F_\chi)=\varnothing$ pour tout~$\chi\in\Ch\Xi_F$, ou bien~$\codim V(F_\chi)=\card F$ pour tout~$\chi\in\Ch\Xi_F$.\hfill$\square$
\bigskip
\begin{coro}
Soit~$f\in{\mathcal S}^*_{n,\mathbb R^n}$, alors~${\mathcal Y}_f=\overline{\re V(f)}={\mathcal F}_f$.
\end{coro}
\pf L'ensemble des spectres d'un SSE constitué par une seule somme d'expo\-nentiel\-les est toujours fermé.\hfill$\square$
\par
\bigskip
Le lemme qui suit concerne le comportement des amibes sous l'action d'un automorphisme~$\mathbb C$-linéaire de~$\mathbb C^n$ qui préserve~$\mathbb R^n$, si~$\varphi$ est un tel automorphisme et~$F$ un SSE à fréquences réelles, on pose
$$
F\circ\varphi:=\{f\circ\varphi\in{\mathcal S}_{n,\mathbb R^n}^*\mid f\in F\}\,.
$$
\begin{lem}\label{cov} Soit~$F$ un SSE à fréquences réelles et~$\varphi:\mathbb C^n\longrightarrow\mathbb C^n$ un
isomorphisme~$\mathbb C$-linéaire tel que~$\varphi(\mathbb R^n)=\mathbb R^n$. Alors~:
\item[$(i)$]
$$
\re \Big[V(F\circ\varphi)\Big]=\re \Big[\varphi^{-1}(V(F))\Big]=\varphi^{-1}(\re V(F))
$$
\item\item[$(ii)$]
$$
{\mathcal Y}_F=\varphi^a({\mathcal Y}_{F\circ\varphi^a})\,,
$$
o\`u~$\varphi^a$ dénote l'adjoint de~$\varphi$ par rapport à la forme hermitienne
standard sur~$\mathbb C^n$.
\end{lem}
{\bf Démonstration.} $(i)$ La première égalité est évidente. Si~$x\in\re\varphi^{-1}(V(F))$ et~$z=x+iy\in \varphi^{-1}(V(F))$, on a~$\varphi(z)=\varphi(x)+i\varphi(y)$ d'o\`u~$\re\varphi(z)=\varphi(x)$,
soit~$x\in\varphi^{-1}(\re V(F))$. D'autre part, si~$x\in\varphi^{-1}(\re V(F))$, il existe un point~$\zeta\in V(F)$ tel que~$\re \zeta=\varphi(x)$ et comme~$\varphi$ est inversible, on a
$$
x=\varphi^{-1}(\re(\zeta))=\re\varphi^{-1}(\zeta)\in\re\varphi^{-1}(V(F))\,.
$$
\par\noindent
$(ii)$~Soit~$f\in F$,~$f(z):=\sum_{\lambda\in\Lambda_f}c_\lambda
e^{\langle z,\lambda\rangle}$, alors, pour tout~$\lambda\in\Lambda_f$ on a
$$
\langle \varphi^a(z),\lambda\rangle
=
\overline{\langle\lambda,\varphi^a(z)\rangle}
=
\overline{\langle\varphi(\lambda),z\rangle}
=
\langle z,\varphi(\lambda)\rangle\,,
$$
d'o\`u
$$
f\circ\varphi^a(z)
=
\sum_{\varphi(\lambda)\in\varphi(\Lambda_f)}
c_{\varphi(\lambda)}
e^{\langle z,\varphi(\lambda)\rangle}
$$
et
$$
\Ch\Xi_{F\circ\varphi^a}
=
\{\chi\circ\varphi^{-1}_{\vert\varphi(\Xi_F)}\mid\chi\in\Ch\Xi_F\}\,.
$$
Ceci fait que, pour tout~$f\in F$ et tout~$\chi\in\Ch\Xi_F$, on ait
$$
f_\chi\circ\varphi^a
=
(f\circ\varphi^a)_{\chi\circ\varphi^{-1}}
$$
ainsi on en déduit 
$$
\re V\big(F_\chi\circ\varphi^a\big)
=
\re V\Big((F\circ\varphi^a)_{\chi\circ\varphi^{-1}}\Big)
$$
et, gr\^ace à~$(i)$
$$
\re V(F_\chi)
=
\varphi^a\Big(\re V\Big((F\circ\varphi^a)_{\chi\circ\varphi^{-1}}\Big)\Big)\,,
$$
d'où la conclusion en prenant l'union sur~$\chi\in\Ch\Xi_F$.\hfill$\square$
\bigskip
\begin{lem}\label{cov2} Soit~$F$ un {\rm SSE}
à fréquences réelles tel que le rang de~$\Xi_F$ soit égale à~$\dim_{\mathbb R}({\rm vect}_{\mathbb R}\Xi_F)$.
Alors on a l'égalité
$$
{\mathcal Y}_F=\re V(F)\,.
$$
\end{lem}
\pf On suppose d'abord que, pour tout~$f\in F$,~$\Lambda_f\subset\mathbb Z^n$ et que~$\Xi_F$ est
de la forme
$$
\Xi_F=\{(m_1,\ldots,m_s,0,\ldots,0)\in\mathbb R^n\mid m_1,\ldots,m_s\in\mathbb Z\}\,,
$$
où~$s\in\{1,\ldots,n\}$ dénote le rang de~$\Xi_F$. Dans ce cas, pour déterminer l'amibe de~$F$ on peut
utiliser les caractères du groupe~$\mathbb Z^n$~;
donc si~$\chi\in\Ch \mathbb Z^n$ est le caractère associé
au~$n$-uplet~
$
(e^{i\theta_1},\ldots,e^{i\theta_n})\,,
$
o\`u~$(\theta_1,\ldots,\theta_s)\in\mathbb R^n$,~$f\in F$ et~$\lambda\in\Lambda_f$,
pour tout~$z\in\mathbb C^n$, on a
\begin{eqnarray*}
\chi(\lambda)e^{\langle z,\lambda\rangle}
&=&
e^{i(m_1\theta_1+\cdots+m_s\theta_s)}e^{z_1 m_1+\cdots+ z_s m_s}\\
&=&
e^{(z_1+i\theta_1)m_1+\cdots+(z_r+i\theta_s)m_s}\\
&=&
e^{\langle z+i\theta,\lambda\rangle}\,,
\end{eqnarray*}
donc~$f_\chi(z)=f(z+i\theta)$. On en tire que, pour tout~$\chi\in\Ch\mathbb Z^n$,~$z\in V(F_\chi)$ si et seulement si~$z+i\theta\in V(F)$, d'o\`u~$\re V(F_\chi)=\re V(F)$ et pour le choix arbitraire de~$\chi\in\Ch \mathbb Z^n$, on déduit que~${\mathcal Y}_F=\re V(F)$.
\par
On passe maintenant au cas général. Supposons que le rang~$s$ de~$\Xi_F$ soit égal à~$\dim_{\mathbb R}
({\rm vect}_{\mathbb R} \Xi_F)$ et soit~$\{\omega_1,\ldots,\omega_s\}$
un système libre de générateurs de~$\Xi_F$.
Les éléments~$\omega_1\ldots,\omega_s$ sont linéairement indé\-pendants sur~$\mathbb R$
car autrement le sous-espace vectoriel de~$\mathbb R^n$ qu'ils engendrent,
à savoir le sous-espace~${\rm vect}_{\mathbb R}\Xi_F$,
aurait dimension plus petite que~$s$. Ceci nous permet de compléter
le système~$\{\omega_1,\ldots,\omega_s\}$ en une
base~$\{\omega_1,\ldots,\omega_s,\omega_{s+1},\ldots,\omega_n\}$ de~$\mathbb R^n$.
Soit~$A$ la matrice donnée par
$$
A:=
\pmatrix{\omega_{11}&\cdots&\omega_{n1}\cr
         \vdots     &\ddots&\vdots  \cr
         \omega_{1n}&\cdots&\omega_{nn}\cr}
\,,
$$
alors, si~$B$ est l'inverse de~$A$ et~$\varphi$ l'automorphisme~$\mathbb C$-linéaire de~$\mathbb C^n$
ré\-pre\-senté dans les bases canoniques par la matrice~$B$,
on voit que~$\varphi(\mathbb R^n)=\mathbb R^n$, et qu'à moins d'une permutation impaire des premières~$s$ colonnes de~$A$ on peut supposer~$\det\varphi>0$. Ceci implique l'égalité
$$
\varphi(\Xi_F)=\{(m_1,\ldots,m_s,0,\ldots,0)\in\mathbb R^n\mid m_1,\ldots,m_s\in\mathbb Z\}\,;
$$
donc, avec les notations du \lemref{cov}, la première partie de la dé\-monstra\-tion nous assure que
$$
{\mathcal Y}_{F\circ\varphi}=\re V(F\circ\varphi)\,,
$$
et un recours au \lemref{cov} nous donne
$$
{\mathcal Y}_F
=
\varphi({\mathcal Y}_{F\circ \varphi})
=
\varphi(\re V(F\circ\varphi))
=
\varphi(\varphi^{-1}(\re V(F)))
=
\re V(F)\,,
$$
ce qui achève la preuve.\hfill$\square$

\bigskip

\begin{coro}\label{cov3}
Si~$F$ est un SSE à fréquences rationnelles alors
$$
{\mathcal Y}_F=\re V(F)\,.
$$
\end{coro}
\pf 
Au vu de \lemref{cov2}, il suffit de vérifier que le rang~$s$ de~$\Xi_F$ est égal à~$
\dim_{\mathbb R}({\rm vect}_{\mathbb R}\Xi_F)\,.
$
Pour cela, soit~$\{\omega_1,\ldots,\omega_s\}\subset\mathbb Q^n$ un
système libre de générateurs de~$\Xi_F$.
On a ~$s\leqslant n\,$; en effet, si~$j\in\{1,\ldots,s\}$ et
$$
\omega_j=\bigg({p_{j1}\over q_{j1}},\ldots,{p_{jn}\over q_{jn}}\bigg)\,,
$$
avec~$p_{j1},\ldots,p_{jn}\in\mathbb Z$ et~$q_{j1},\ldots,q_{jn}\in\mathbb Z^*$,
alors, en posant
$$
\mu:={\textsc{ppmc}}\{q_{jk}\in\mathbb Z\mid j\in\{1,\ldots,s\}\;,k\in\{1,\ldots,n\}\}\,;
$$
on voit que~$\mu\neq 0$, donc~$\Xi_f$ est isomorphe à~$\mu\Xi_f$ et
comme~$\mu\Xi_F\subseteq\mathbb Z^n$, on en tire que~$s\leqslant n$.
De plus,~$\omega_1,\ldots,\omega_s$ sont~$\mathbb Q$-linéairement indépendants
car en multipliant une éventuelle relation de dépendance linéaire sur~$\mathbb Q$ par
le plus petit multiple commun des dénominateurs des coefficients de la relation, on obtient
une relation sur~$\mathbb Z$, ce qui est contraire au fait que les éléments~$\omega_1,\ldots,\omega_s$
définissent une famille libre sur~$\mathbb Z$. Par conséquent, on peut
compléter~$\{\omega_1,\ldots,\omega_s\}$ en une
base~$\{\omega_1,\ldots,\omega_s,\omega_{s+1},\ldots,\omega_n\}$ de~$\mathbb Q^n$.
Comme dans la démonstration du \lemref{cov2}, soit~$A$ la matrice donnée par
$$
A:=
\pmatrix{\omega_{11}&\cdots&\omega_{n1}\cr
         \vdots     &\ddots&\vdots  \cr
         \omega_{1n}&\cdots&\omega_{nn}\cr}
\,,
$$
alors, si~$B$ est l'inverse de~$A$ et~$\varphi$ l'automorphisme~$\mathbb C$-linéaire de~$\mathbb C^n$
re\-pré\-senté dans les bases canoniques par la matrice~$B$,
on a que~$\varphi(\mathbb Q^n)=\mathbb Q^n$ et, à moins d'une permutation impaire des premières~$s$ colonnes de~$A$, on peut supposer~$\det\varphi>0$. Ceci implique l'égalité
$$
\varphi({\rm vect}_{\mathbb Q}\Xi_F)=\{(m_1,\ldots,m_s,0,\ldots,0)\in\mathbb R^n\mid m_1,\ldots,m_s\in\mathbb Q\}\,,
$$
o\`u~${\rm vect}_{\mathbb Q}\Xi_F$ dénote le~$\mathbb Q$-sous-espace vectoriel
de~$\mathbb Q^n$ engendré par~$\Xi_F$, d'o\`u 
$$
\dim_{\mathbb R}({\rm vect}_{\mathbb R}\Xi_F)
=
\dim_{\mathbb R}(\varphi({\rm vect}_{\mathbb R}\Xi_F))
=
\dim_{\mathbb Q}(\varphi({\rm vect}_{\mathbb Q}\Xi_F))
=
s\,,
$$
ce qui conclut la preuve.\hfill$\square$
\bigskip
\begin{rmq} 
{\rm
Le \corref{cov3} a comme conséquence le fait que notre notion d'amibe pour un
système de sommes d'exponentielles généralise la notion classique d'amibe.
En fait si~$P=\{p_1,\ldots,p_r\}\subset\mathbb C[x_1^{\pm 1},\ldots,x_n^{\pm 1}]$ est
un système de polynômes de Laurent non nuls,~$V(P)$ son ensemble
des zéros dans le tore~$(\mathbb C^*)^n$ et~${\mathcal A}_P:={\rm Ln} V(P)$
son amibe au sens classique, la substitution~$x_j=e^{z_j}$, pour~$j=1,\ldots,n$,
transforme~$P$ en le SSE à fréquences entières~$F:=\{f_1,\ldots,f_r\}
\subset{\mathcal S}_{n,\mathbb Z^n}^*$, o\`u, pour~$1\leqslant k\leqslant r$ et~$z\in\mathbb C^n$, on pose
$$
f_k(z):=p_k(e^{z_1},\ldots,e^{z_n})\,.
$$
Comme, pour tout~$j=1,\ldots,n$,
~$
\ln\vert x_j\vert=\ln\vert e^{z_j}\vert=\re z_j\,
$,
on en déduit que
$$
{\mathcal A}_P={\mathcal Y}_F={\mathcal F}_F\,.
$$}
\end{rmq}
\bigskip
\begin{exe}\label{bal1}
{\rm 
Soit~$\gamma\in\mathbb R\setminus\mathbb Q$ et~$f\in{\mathcal S}_{1,\mathbb R}^*$ donnée, pour~$z\in\mathbb C$, par
\begin{eqnarray*}
f(z)
&=&
\cos(iz)+\sin(i\gamma z)-2\\
&=&
{1\over 2}\big(e^{-z}+e^z\big)+{1\over 2i}\big(e^{-\gamma z}-e^{\gamma z}\big)-2\,.
\end{eqnarray*}
L'ensemble~$\re V(f)$ n'est pas fermé dans~$\mathbb R$ donc~$\re V(f)\subsetneq{\mathcal Y}_f$. En effet,
si~$z$ est imaginaire pur,~$f(z)=0$ si et seulement
si~$\cos iz=1$ et~$\sin (i\gamma z)=1$, soit si et seulement si
$$
iz\in
2\,\pi\mathbb Z
\cap
\Big((\pi/ 2\gamma)+(2\,\pi/\gamma)\mathbb Z\Big)=\varnothing\,,
$$
en particulier~$0\notin\re V(f)$. D'autre part, si~$\chi\in\Ch\Xi_f$ est tel que~$\chi(1)=1$ et~$\chi(\gamma)=-i$, on a bien
$$
f_\chi(z)=\cos(iz)+\cos(i\gamma z)-2\,,
$$
d'où~$f_\chi(0)=0$ et donc~$0\in{\mathcal Y}_f=\overline{\re V(f)}$.\footnote{Une preuve directe de ceci n'utilisant pas le langage des amibes m'a été signalée par Michel Balazard.} Puisque le rang du groupe~$\Xi_f$ est égale à~$2$ on voit que le \lemref{cov2} est en général faux si le rang de~$\Xi_F$ est plus grand que~$\dim_{\mathbb R}({\rm vect}_{\mathbb R}\,\Xi_F)$.}
\end{exe}
\begin{exe}\label{jam1}
{\rm 
Soit~$\gamma\in\mathbb R\setminus\mathbb Q$ et~$f\in{\mathcal S}^*_{1,\mathbb R}$ donnée, pour~$z\in\mathbb C$, par
$$
f(z)
=
(e^z-1)(e^{\gamma z}-e^\gamma)\,,
$$
alors~$\re V(f)=\{0,1\}={\mathcal Y}_f$ malgré les hypothèses du~\lemref{cov2} ne soient pas satisfaites. La condition énoncée dans le~\lemref{cov2} est donc suffisante mais pas nécessaire pour qu'on ait~${\mathcal Y}_f=\re V(f)$.} 
\end{exe}
\begin{exe}\label{jam2}
{\rm Soit~$\gamma\in\mathbb R\setminus\mathbb Q$ et~$F=\{f,g\}\subset{\mathcal S}_{1,\mathbb R}^*\,$, où~$f$ et~$g$ sont données, pour~$z\in\mathbb C$, par
$$
f(z)
=
\cos(iz)+\sin(i\gamma z)-2
\qquad{\rm et}\qquad
g(z)
=
e^{z}-1\,.
$$
Il s'agit d'un système qui n'a pas de solutions (car~$f$ n'a pas de zéros imaginaires purs alors que~$g$ n'a que de tels zéros), donc~$\overline{\re V(F)}=\varnothing$. D'autre part~${\mathcal Y}_F\neq\varnothing$ car~$0\in V(F_\chi)$, où~$\chi$ désigne le caractère tel que~$\chi(1)=1$ et~$\chi(\gamma)=-i$. Donc il existe bien de systèmes~$F\subset{\mathcal S}^*_{n,\mathbb R^n}$ qui n'ont pas de zéros et dont l'amibe~${\mathcal Y}_F$ n'est pas vide. Dans ces cas l'amibe~${\mathcal Y}_F$ est trop grande (donc peu intéressante) et le~\thmref{proprF} est faux.}
\end{exe}

\section{k-convexité selon Henriques.}
\par
Dans cette section on va faire quelques remarques autour de la notion de~$k$-convexité pour un ouvert d'un espace affine réel telle qu'elle a été introduite dans \cite{Hen}, auquel on renvoie pour toutes les définitions, les détails techniques et tous les résultats que l'on évoquera dans la suite, en particulier en ce qui concerne le complexe des chaînes polyédrales. 
\par
Si~$\varnothing\neq X\subset\mathbb R^n$ est un ouvert, on note~${}^{\rm pl}C_\bullet(X)$ le complexe des chaînes poly\-édrales de~$X$, il est obtenu comme le quotient du complexe~${}^\Delta C_\bullet(X)$ de chaînes linéaires par morceaux de~$X$ modulo la relation~$\sim$ d'équivalence géométrique de ces chaînes. Si~$\sigma=\sum_{j=1}^m\lambda_j\sigma_j\in{}^{\rm pl}C_k(X)$, avec~$\lambda_j\neq 0$ pour tout~$j$ et~$c=[\sigma]_\sim\in{}^\Delta C_k(X)$, on rappelle que le support~$\Supp\sigma$ de~$\sigma$ est l'union des images des chaînes~$\sigma_j$ qui apparaissent dans l'expression de~$\sigma$ et que
$$
\Supp c:=\bigcap_{\tau\sim\sigma} \Supp\tau\,,
$$
ce dernier étant bien défini en vertu du Lemme~2.4 dans \cite{Hen}. On rappelle aussi que l'homologie du complexe~${}^\Delta C_\bullet(X)$ est isomorphe à l'homologie singulière de~$X$, (\cite{Hen} Lemme~2.2), donc dans toute question de~$k$-convexité pour un ouvert~$X$ d'un espace affine réel, on pourra utiliser l'homologie du complexe~${}^\Delta C_\bullet(X)$ au lieu de celui des chaînes singulières de~$X$.
\par
Le terme~$k$-{\sl convexité} n'est pas nouveau en Mathématiques, il existe en fait en analyse complexe de plusieurs variables ainsi qu'en analyse fonctionnelle. Néanmoins ces notions analytiques ne ressemblent pas à la notion présentée par Henriques, qui me parait quand même assez nouvelle. On mentionne d'ailleurs que Mikalkhin~\cite{Mi2} a introduit, sous le même nom de~$k$-convexité, une notion plus forte que celle d'Henriques.
\par
Il faut remarquer que, si~$k\in\mathbb N$ est fixé, la~$k$-convexité dans~$\mathbb R^n$ ne dévient intéressante que pour~$n\geqslant k+2$, sinon tout sous-ensemble de~$\mathbb R^n$ est~$k$-convexe. Des simples exemples sont le complémentaire d'une union finie de droites dans~$\mathbb R^3$, qui est~$1$-convexe mais qu'il n'est pas~$0$-convexe\footnote{Un autre exemple assez explicatif d'un tel sous-ensemble m'a été signalé par Mikael Passare, il s'agit du complémentaire d'une ``tour Eiffel'' dans~$\mathbb R^3$.} et, plus en général, le complémentaire d'une union finie de~$k$-sous-espaces affines dans~$\mathbb R^{k+2}$, qui est~$k$-convexe mais il n'est pas~$\ell$-convexe, pour~$\ell<k$. Par contre, le complémentaire d'un ensemble fini de points dans~$\mathbb R^3$ n'est pas~$0$-convexe, ni~$1$-convexe (mais il est trivialement~$2$-convexe).
\par
La ``faiblesse'' de la notion de~$k$-convexité croit avec~$k$.
\begin{lem}
Soit~$X\subset\mathbb R^n$ un sous-ensemble non vide et soit~$k\in\mathbb N$. Alors, si~$X$ est~$k$-convexe, il est aussi~$(k+1)$-convexe.
\end{lem}
\pf On suppose par l'absurde que~$X$ soit~$k$-convexe mais qu'il ne soit pas~$(k+1)$-convexe. Il existe donc un~$(k+2)$-sous-espace affine orienté~$S$ de~$\mathbb R^n$ qui rencontre~$X$ et il existe aussi une classe non nulle~$c$ dans~$\widetilde H^+_{k+1}(S\cap X,\mathbb Z)$ dont l'image, (sous le morphisme induit par l'inclusion), dans~$\widetilde H_{k+1}(X,\mathbb Z)$ est nulle. Soit alors~$\sigma$ un~$(k+1)$-cycle non négatif dans~$S\cap X$ qui représente la classe~$c$ et~$S^\prime$ est un~$(k+1)$-sous espace affine orienté de~$S$ tel que
l'intersection~$\sigma^\prime:=S^\prime\cap \sigma$ soit un~$k$-cycle non négatif et non nul contenu dans~$S^\prime\cap X$, (un tel sous-espace existe car autrement~$\sigma$ représenterait la classe nulle de~$\widetilde H^+_{k+1}(S\cap X,\mathbb Z)$). Puisque~$\sigma$ représente la classe nulle dans~$\widetilde H_{k+1}(X,\mathbb Z)$, on déduit que~$\sigma^\prime$ représente la classe nulle dans~$\widetilde H_{k}(X,\mathbb Z)$, ce qui est contraire à la~$k$-convexité de~$X$.\hfill$\square$ 

\bigskip
On termine la section par le lemme suivant.
\begin{lem} Soit~$\varphi:\mathbb R^n\longrightarrow\mathbb R^n$
un isomorphisme d'espaces affines qui préserve l'orientation.
Alors si~$X\subset \mathbb R^n$ est~$k$-convexe,~$\varphi(X)$ l'est.
\end{lem}
\pf Si~$X=\varnothing$ il n'y a rien à montrer. Sinon, la restriction de~$\varphi$ à~$X$ induit un homéomorphisme de~$X$ sur~$\varphi(X)$ donc un isomorphisme en homologie réduite~$\varphi_*:{\tilde H}_\bullet(X)
\longrightarrow{\tilde H}_\bullet(\varphi(X))$. En outre, comme~$\varphi$ préserve
l'orientation, pour tout~$(k+1)$-sous-espace affine orienté~$S$ de~$\mathbb R^n$ qui rencontre~$X$,~$\varphi(S)$ est un sous-espace affine de~$\mathbb R^n$ qui est
isomorphe à~$S$, en tant qu'espace affine réel orienté, et qui rencontre~$\varphi(X)$~; d'autre part, tout~$(k+1)$-sous-espace affine orienté
de~$\mathbb R^n$ qui rencontre~$\varphi(X)$ est de la forme~$\varphi(S)$
pour un unique~$S$. Enfin, pour tout~$(k+1)$-sous-espace
affine~$S$ de~$\mathbb R^n$ qui rencontre~$X$ et tout~$x\in S\setminus X$, l'isomorphisme~$\varphi$ induit un isomorphisme
$$
\varphi_*:\mathbb Z={\tilde H}_k(S\setminus\{x\})\longrightarrow
{\tilde H}_k(\varphi(S)\setminus\{\varphi(x)\})=\mathbb Z
$$
qui, comme on le voit facilement, n'est rien d'autre que l'isomorphisme identité. On peut donc conclure
la démonstration, en fait, pour tout~$(k+1)$-sous-espace affine orienté~$S$ de~$\mathbb R^n$ qui rencontre~$X$ et tout~$x\in S\setminus X$,
$$
{\tilde H}_k^+(\varphi(S)\cap\varphi(X))\setminus\{0\}
=\varphi_*({\tilde H}^+_k(S\cap X)\setminus\{0\})
$$
et de plus le diagramme suivant
$$
\matrix{
{\widetilde H}_k(\varphi(S)\cap \varphi(X))
&
\hfl{}{}
&
{\widetilde H}_k(\varphi(S)\setminus\{\varphi(x)\})\cr
&&&\cr
\vfl{\varphi_*^{-1}}{}
&
&
\vfl{}{id}\cr
&&&\cr
{\widetilde H}_k(S\cap X)
&
\hfl{}{}
&
{\widetilde H}_k(S\setminus\{x\})=\mathbb Z\cr
}
$$
(où les flèches horizontales sont induites par l'inclusion), est commutatif.\hfill$\square$

\section{Le complémentaire de l'amibe.}

Dans cette section on démontre un résultat sur le
complé\-mentaire~${\mathcal F}_F^c$ de l'amibe d'un SSE~$F$ à fréquences réelles
qui constitue le pendant du \thmref{hen}.
Pour cela, on aura besoin d'une hypothèse géométrique sur les fréquences de~$F$, à savoir
l'hypothèse que l'ensemble des spectres de~$F$ soit fermés.
\par
\medskip
\begin{thm}
Soit~$F\subset{\mathcal S}_{n,\mathbb R^n}^*$ un SSE dont l'ensemble des spectres est fermé. Si~$F$ est constitué par~$(k+1)$ sommes d'exponentielles, le complé\-mentaire~${\mathcal F}_F^c$ de l'amibe de~$F$
est un sous-ensemble~$k$-convexe de~$\mathbb R^n$.
\end{thm}
\pf 
L'ensemble des spectres de~$F$ est fermé donc~${\mathcal F}_F={\mathcal Y}_F$ et,
pour tout~$\chi\in\Ch \Xi_F$, le SSE~$F_\chi$ est régulier.
Si~$\dim_{\mathbb C}(\aff_{\mathbb C}\Gamma_F)<(k+1)$, pour tout~$\chi\in\Ch\Xi_F$,
on a~$V(F_\chi)=\varnothing\,$, donc~${\mathcal Y}_F^c=\mathbb R^n$ qui est évidemment~$k$-convexe.
Par contre, si~$\dim_{\mathbb C}(\aff_{\mathbb C}\Gamma_F)\geqslant (k+1)$, l'ensemble
analytique~$V(F_\chi)$ est non vide et de codimension~$(k+1)$ dans~$\mathbb C^n$,
pour tout~$\chi\in\Ch\Xi_F$. On conduit la démonstration en trois étapes.
\par
$(i)$ Si~$\Lambda_f\subset\mathbb Z^n$ pour tout~$f\in F$, l'amibe~${\mathcal Y}_F$
co{\"\i}ncide avec l'amibe (au sens classique)~${\mathcal A}_P$
d'un système~$P$ de polynômes de Laurent de~$n$ variables tel que
la codimension, dans~$(\mathbb C^*)^n$, de l'ensemble algébrique~$V(P)$ soit
égale à~$(k+1)$. Gr\^ace au \thmref{hen}, on peut conclure
que~${\mathcal Y}_F^c$ est~$k$-convexe dans ce cas.
\par
$(ii)$ On suppose maintenant que~$\Lambda_f\subset\mathbb Q^n$
pour tout~$f\in F$, et, comme dans la démons\-tration du \corref{cov2},
on peut trouver un automorphisme $\mathbb C$-linéaire~$\varphi$ de~$\mathbb C^n$ tel que~$\det\varphi>0$,~$\varphi(\mathbb R^n)=\mathbb R^n$ et~$\varphi(\Xi_F)\subset\mathbb Z^n$.
Ainsi, avec les mêmes notations qu'au \lemref{cov}, on a
$$
{\mathcal Y}_F=\varphi^a({\mathcal Y}_{F\circ\varphi^a})
\qquad
{\rm et}
\qquad
{\mathcal Y}_F^c=\varphi^a({\mathcal Y}_{F\circ\varphi^a}^c)
$$
car l'adjoint~$\varphi^a$ de~$\varphi$ est aussi bijectif.
En outre, le fait que~$\varphi$ soit un isomorphisme implique que l'ensemble des spectres du système~$F\circ\varphi^a$ soit aussi fermé, donc~$\dim_{\mathbb C}(\aff_{\mathbb C}\Gamma_{F\circ\varphi^a})\ge (k+1)$,
et~$\codim V(F\circ\varphi^a)=(k+1)$.
\par
Or, comme ~$\Xi_{F\circ\varphi^a}=\varphi(\Xi_F)\subset\mathbb Z^n$, la première
partie de la démonstration montre que l'ensemble~${\mathcal Y}^c_{F\circ\varphi^a}$
est~$k$-convexe dans~$\mathbb R^n$ et vu que~$\det\varphi^a>0$ un recours au \lemref{cov} permet de conclure la démonstration dans ce deuxième cas.
\par
$(iii)$ On passe donc au cas général où~$\Lambda_f\subset\mathbb R^n$ 
pour tout~$f\in F$. Si~$\{\omega_1,\ldots,\omega_r\}$ est un système libre de 
générateurs de~$\Xi_F$ on aura
$$
f(z)=\sum_{k\in A_f}a_{f,k} e^{k_1\langle z,\omega_1\rangle+\cdots+k_r\langle z,\omega_r\rangle}\,,
$$
où~$A_f\subset\mathbb Z^r$ est un sous-ensemble fini et~$a_{f,k}\in\mathbb C^*$ pour tout~$k\in A_f$. 
Pour tout~$j\in\{1,\ldots,r\}$, soit~$(\omega_{j,\ell})_{\ell\in\mathbb N}\subset\mathbb Q^n$ une suite convergeante  
vers~$\omega_j$ et, pour tout~$\ell\in\mathbb N$, soit~$F^{[\ell]}:=\{f^{[\ell]}\in{\cal S}_{n,\mathbb R}
\mid f\in F\}$, o\`u~$f^{[\ell]}$ est la somme d'exponentielles donnée par
$$
f^{[\ell]}(z):=
\sum_{k\in A_f} 
a_{f,k} e^{k_1\langle z,\omega_{1,\ell}\rangle+\cdots+k_r\langle z,\omega_{r,\ell}\rangle}\,.
$$
On voit ainsi que, pour tout~$f\in F$, 
la suite des polytopes~$(\Gamma_{f^{[\ell]}})_{\ell\in\mathbb N}$ converge vers le polytope~$\Gamma_f$ 
pour la métrique de Hausdorff; 
par conséquent, pour~$\ell$ assez grand, l'ensemble des spectres du système~$F^{[\ell]}$ est aussi fermé (donc~$F^{[\ell]}$ est régulier) et~$\dim_{\mathbb C}(\aff_{\mathbb C}\Gamma_{F^{[\ell]}})\geqslant (k+1)$. 
Ceci implique que, pour~$\ell$ assez grand, l'ensemble analytique~$V(F^{[\ell]})$ 
est non vide et de codimension~$(k+1)$ dans~$\mathbb C^n$. D'autre part, pour tout~$f\in F$, 
le support de~$f^{[\ell]}$ est contenu dans~$\mathbb Q^n$,
donc, en vertu de la deuxième partie de la démonstration, on sait que 
pour~$\ell$ assez grand, l'ensemble~${\cal Y}_{F^{[\ell]}}^c$ est~$k$-convexe. 
\par
De fa\c con analogue, pour tout~$\chi\in\Ch\Xi_F$ et tout~$\ell\in\mathbb N$, 
on peut définir~$(F_\chi)^{[\ell]}$, et puisque, pour tout~$\chi\in\Ch\Xi_F$, tout~$\ell\in\mathbb N$ et 
tout~$f\in F$, on a~
$
\Lambda_{f^{[\ell]}}=\Lambda_{(f_\chi)^{[\ell]}}\subset\mathbb Q^n
$, 
on peut également conclure que, pour tout caractère~$\chi\in\Ch\Xi_F$ et pour~$\ell$ assez grand,  
l'ensemble~${\cal Y}_{(F_\chi)^{[\ell]}}^c$ est~$k$-convexe.
\par 
Les hypothèses~$F\neq\{0\}$ et~${\mathcal Y}_F\neq\varnothing$ impliquent l'existence d'un~$(k+1)$-sous-espace affine orient\'e~$S$ de~$\mathbb R^n$ tel que~$\varnothing\neq S\cap{\cal Y}_F^c\neq S$. 
Soit $S$ un tel sous-espace affine (d'espace vectoriel sous-jacent $E_S$) et supposons, par l'absurde, qu'il existe une classe~$\gamma\in{\widetilde H}^+_k(S\cap{\cal Y}^c_F)\setminus\{0\}$ dont l'image est nulle sous le morphisme
$$
\iota:{\widetilde H}_k(S\cap{\cal Y}^c_F)\longrightarrow\widetilde H_k({\cal Y}^c_F)\,
$$
induit par l'inclusion; il s'agit de montrer que l'existence d'un tel élément conduit \`a une contradiction. 
On choisit pour cela un représentant~$c$ de~$\gamma$ dans le groupe~${\cal C}^\Delta_k( S\cap {\cal Y}^c_F)$ (c'est-\`a-dire un $k$-cycle  
affine par morceaux de l'ouvert $S\cap {\cal Y}^c_F$ de l'espace affine~$(k+1)$-dimensionnel~$S$); grâce au Lemme~2.7 de~\cite{Hen}, il existe une unique  
$(k+1)$-chaîne affine par morceaux $C$ de ${\cal C}^\Delta_{k+1} (S)$ (dépendant de $c$) 
telle que $\partial C=c$ et l'hypothèse que la classe d'homologie de $c$ dans~$S\cap {\cal Y}_F^c$ soit 
non nulle équivaut (pour le même Lemme~2.7 de~\cite{Hen}) à ce que le support de $C$ ne soit pas inclus dans~${\cal Y}_F^c$; il 
existe donc un caractère $\chi_o$ de $\Xi_F$ tel que le support de $C$ n'est pas inclus dans~$S\cap({\rm Re}\, V(F_{\chi_o}))^c$.
\\
En outre, comme~$\Supp c\subset{\cal Y}^c_F$, on voit que la classe nulle de~$\widetilde H_k({\cal Y}^c_F)$ peut être représentée par le cycle~$c$, donc il existe un élément~
$D\in{\cal C}^\Delta_{k+1}({\cal Y}^c_F)$, tel que~$\partial D=c$ dans~${\cal Y}^c_F$.
\par
On admet pour l'instant qu'il existe~$L\in\mathbb N$ tel que pour tout~$\ell\geqslant L$ on a 
$$
\Supp c \cup\Supp D\subseteq {\cal Y}_{(F_{\chi})^{[\ell]}}^c\,,
\eqno{(*)_\chi^\ell}
$$
pour tout~$\chi\in\Ch\Xi_F$, donc, pour~$\ell\geqslant L$, la relation~$(*)_{\chi_o}^\ell$ implique que~$c$ représente une classe d'homologie~$\gamma_{\chi_o,\ell}$ de~$\widetilde H_k(S\cap {\cal Y}^c_{(F_{\chi_o})^{[\ell]}})$ dont l'image est nulle sous le morphisme
$$
\iota_\ell:
\widetilde H_k(S\cap {\cal Y}^c_{(F_{\chi_o})^{[\ell]}})
\longrightarrow
\widetilde H_k({\cal Y}^c_{(F_{\chi_o})^{[\ell]}})\;,
$$
induit par l'inclusion. En outre, l'hypothèse~$\gamma\in{\widetilde H}^+_k(S\cap{\cal Y}^c_F)$ implique que, si~$\ell\geqslant L$, on a~$\gamma_{\chi_o,\ell}\in
\widetilde H^+_k(S\cap {\cal Y}^c_{(F_{\chi_o})^{[\ell]}})$. En fait si, pour~$\ell\geqslant L$ et~$x$ appartenant~$S\setminus{\cal Y}^c_{(F_{\chi_o})^{[\ell]}}$,~$\upsilon_x$ dénote le générateur standard du groupe de cohomologie de de~Rham~$H^k_{dR}(S\setminus\{x\})$, 
$$
\upsilon_x:=
{1\over \varkappa_k}
\sum_{j=0}^k
(-1)^j{\xi_j-x_j\over\parallel \xi-x\parallel^{k+1}}
\,d\xi_{[j]}\,,
$$
($\varkappa_k$ étant le volume~$k$-dimensionnel de la sphère~$k$-dimensionnelle), on a
$$
\int_c \upsilon_x >0\quad{\rm lorsque}~x\in\Supp C\,,
$$
ou alors
$$
\int_c \upsilon_x=0\quad
{\rm lorsque}~x\notin\Supp C\,.
$$
Si l'on change de représentant pour~$\gamma_{\chi_o,\ell}$, il est facile de voir (par le théorème de Stokes) qu'aucune des deux intégrales ci-dessus peut devenir négative, donc, grâce au Lemme~3.2 de~\cite{Hen}, si~$\ell\geqslant L$, la classe~$\gamma_{\chi_o,\ell}$ est non négative dans~$\widetilde H_k(S\cap {\cal Y}^c_{(F_{\chi_o})^{[\ell]}})$, d'autre part, si~$\ell\geqslant L$, la~$k$-convexité de~${\cal Y}^c_{(F_{\chi_o})^{[\ell]}}$ implique que~$\gamma_{\chi_o,\ell}$ représente la classe nulle dans le groupe~${\widetilde H}_k(S\cap {\cal Y}^c_{(F_{\chi_o})^{[\ell]}})$, soit~$\Supp C\subseteq{\cal Y}^c_{(F_{\chi_o})^{[\ell]}}$.
\par 
La contradiction attendue viendra alors du fait que l'on sait que le support de $C$ n'est pas inclus dans $S\cap {\rm Re}\, (V(F_{\chi_o}))^c$. En fait on peut trouver un point~$x\in S\cap {\rm Re}\, (V(F_{\chi_o}))$ qui appartient aussi \`a l'intérieure relatif de~$\Supp C$, donc il existe un voisinage~$W$ de~$x$ tel que~$W\cap S$ soit entièrement contenu dans~$\Supp C$. Si~$y\in\mathbb R^n$ est tel que~$x+iy\in V(F_{\chi_o})$, l'intersection de~$V(F_{\chi_o})$ avec
$$
U:=S+i(y+E_S)\,,
$$
constitue un ensemble analytique discret dans~$\mathbb C^n$. Soit donc~$B$ dans~$\mathbb C^n$ une boule ouverte de centre~$x+iy$ qui ne contient pas d'autres points de~$V(F_{\chi_o})\cap U$. Pour~$\ell$ assez grand, l'ensemble analytique~$V((F_{\chi_o})^{[\ell]})\cap U$ est aussi discret et, dans ce cas, la version en plusieurs variables du théorème de Rouché assure que cet ensemble admet dans~$B$ le même nombre d'éléments que~$V(F_{\chi_o})\cap U$ y admet, soit un seul élément, que l'on note~$x_\ell+iy_\ell$. 
\par
Il est clair que la suite des points~$x_\ell+iy_\ell$ tend vers~$x+iy$ et, en particulier, que les points de la suite~$(x_\ell)$ appartiennent à~$W\cap S$, pour~$\ell$ assez grand. Mais alors on a trouvé la contradiction attendue, car, pour~$\ell$ assez grand, on a d'une part~$x_\ell\in{\cal Y}_{(F_{\chi_o})^{[\ell]}}$ et d'autre part~$x_\ell\in W\cap S\subset\Supp C\subset {\cal Y}_{(F_{\chi_o})^{[\ell]}}^c$.
\par
Pour terminer la démonstration, il nous reste \`a prouver qu'il existe~$L\in\mathbb N$ tel que, pour tout~$\ell\geqslant L$, la relation~$(*)_\chi^\ell$ 
est vérifiée pour tout~$\chi\in\Ch\Xi_F$. On commence par 
remarquer qu'il existe un nombre fini~$m$ de boules fermées $\overline B(x_s,\varepsilon_{x_s})$,~$1\leqslant s\leqslant m$, telles que 
$$
{\rm Supp}\, c \cup \Supp D\subset \bigcup\limits_{s=1}^m \overline B(x_s,\varepsilon_{x_s})  
\subset {\cal Y}^c_F\,. 
$$
Il suffit de montrer que, pour chaque $1\leqslant s\leqslant m$, il existe 
un~$l_s\in\mathbb N$ tel que, pour tout entier~$\ell\geqslant l_s$, on ait 
$$
\overline B(x_s,\varepsilon_{x_s})\subset{\cal Y}_{(F_\chi)^{[\ell]}}^c\,,
$$
pour tout $\chi \in \Xi_F$, et prendre en suite~$L:=\max\{l_s\mid 1\leqslant s\leqslant m\}$. 
On prouve ceci par l'absurde$\,$; on suppose que pour un certain~$s$,~$1\leqslant s\leqslant m$, il existe un une suite strictement croissante $(\ell_q)\subseteq\mathbb N$ et une suite~$(\chi_q)\subset\Ch\Xi_F$ 
telles que 
$$
\overline B(x_s,\varepsilon_{x_s})\cap{\cal Y}_{(F_{\chi_q})^{[\ell_q]}}\neq\varnothing\,. 
$$  
Comme, pour tout~$q\in\mathbb N$,~${\cal Y}_{(F_{\chi_q})^{[l_q]}}=\re V((F_{\chi_q})^{[l_q]})\,$, 
on déduit l'existence d'une suite de points~$\xi_q$ de 
$\overline B(x_s,\varepsilon_{x_s})$ et d'une suite de points~$\eta_q$ de $\mathbb R^n$ tels que, 
pour tout~$f\in F$ et tout~$q\in\mathbb N$, on ait
$$
(f_{\chi_q})^{[\ell_q]}(\xi_q+i\eta_q)=0\,,
$$
soit
$$
(f_{\tilde \chi_q})^{[\ell_q]}(\xi_q)=0\,,
$$
o\`u, pour tout~$q\in\mathbb N\,$,~$\tilde \chi_q:=\chi_q \kappa_q$, 
~$\kappa_q$ désignant le caractère de~$\Xi_F$ donné, pour~$1\leqslant j\leqslant r$, par~
$
\kappa_q(\omega_j)=e^{i\langle\eta_q,\omega_{j,\ell_q}\rangle}\,. 
$
Par compacité de~$\overline B(x_s,\varepsilon_{x_s})$ et 
de~$\Ch\Xi_F$, on extrait une sous-suite~$(\xi_{q_r})$ 
et une sous-suite~$(\tilde \chi_{q_r})$ convergeantes respectivement 
vers un point~$\tilde\xi$ de la boule~$\overline B(x_s,\varepsilon_{x_s})$ 
et un caractère~$\tilde \chi$ de~$\Ch\Xi_F$; en passant à la limite, on a 
donc, pour tout~$f\in F$, 
$$
0
=
\lim_{r\to\infty}(f_{\tilde \chi_{q_r}})^{[\ell_{q_r}]}(\xi_{q_r})
=
f_{\tilde \chi}(\tilde \xi)\,,
$$
ce qui est absurde, vu que~$\tilde \xi \in \overline B(x_s,\varepsilon_{x_s})\subset{\cal Y}_F^c$.\hfill$\square$

\bigskip
\par

Je remercie mon directeur de thèse, Alain Yger, pour le support qu'il m'a témoigné pendant la préparation de cet article, ainsi que Michel Balazard et Mikael Passare pour les exemples fournis.

\vfill\eject

\bigskip
\noindent
\small
James SILIPO
\\
LaBAG, Institut de Mathématiques
\\
U.F.R. de Mathématiques et Informatique, Université Bordeaux 1
\\
351 cours de la Libération, 33405, Talence Cedex, France
\\
{\it silipo@math.u-bordeaux1.fr}

\vfill\eject
\end{document}